\documentclass[fleqn]{mat01}
\usepackage{times,mathtimy,amssymb}
\begin{document}

\setcounter{page}{365}
\firstpage{365}

\font\xx=msam5 at 10pt
\def\ab{\mbox{\xx{\char'03}}}

\newtheorem{theo}{Theorem}
\renewcommand\thetheo{\arabic{section}.\arabic{theo}}
\newtheorem{theor}[theo]{\bf Theorem}
\newtheorem{lem}[theo]{Lemma}
\newtheorem{propo}[theo]{\rm PROPOSITION}
\newtheorem{rema}[theo]{Remark}
\newtheorem{defn}[theo]{\rm DEFINITION}
\newtheorem{exam}{Example}
\newtheorem{pol}[theo]{Proof of Lemma}
\newtheorem{coro}[theo]{\rm COROLLARY}
\newtheorem{claim}[theo]{Claim}
\newtheorem{conjecture}[theo]{Conjecture}

\renewcommand{\theequation}{\thesection\arabic{equation}}

\newcommand{\trd}{\rm Trd}
\newcommand{\nrd}{\rm Nrd}
\newcommand{\id}{\rm Id}
\newcommand{\Int}{\rm Int}
\newcommand{\End}{\rm End}
\newcommand{\br}{\rm Br}
\newcommand{\Sym}{\rm Sym}
\newcommand{\Skew}{\rm Skew}
\newcommand{\GL}{\rm GL}
\newcommand{\Spin}{\rm Spin}
\newcommand{\Gal}{\rm Gal}

\renewcommand{\H}{\mathbb{H}}

\newcommand{\R}{\mathbb{R}}
\newcommand{\wt}[1]{\wedge^2 #1}
\newcommand{\iso}{\xrightarrow{\sim}}
\newcommand{\ad}{{\mathrm {ad}}}
\newcommand{\ga}{\gamma}
\newcommand{\s}{\sigma}
\newcommand{\la}{\lambda}
\newcommand{\ot}{\otimes}
\newcommand{\mz}{{\mathbb Z}}
\newcommand{\ra}{\rightarrow}
\newcommand{\sig}{\sigma}
\newcommand{\sq}[1]{#1^\star/#1^{\star 2}}
\newcommand{\ta}{\theta}
\newcommand{\w}{\wedge}
\newcommand{\leven}{\la^{\mathrm{even}}}
\newcommand{\IF}[1]{{I^{#1}\negthinspace F}}
\newcommand{\qform}[1]{{\langle{#1}\rangle}}
\newcommand{\pform}[1]{{\langle\langle{#1}\rangle\rangle}}

\title{Pfister involutions}

\markboth{E~Bayer-Fluckiger, R~Parimala and A~Qu\'eguiner-Mathieu}{Pfister involutions}

\author{E~BAYER-FLUCKIGER$^{*}$, R~PARIMALA$^{\dagger}$ and A~QU\'EGUINER-MATHIEU$^{\ddagger}$}

\address{$^{*}$D\'epartement de Math\'ematiques,
Ecole Polytechnique, F\'ed\'erale de Lausanne, 1015 Lausanne, Switzerland\\
\noindent $^{\dagger}$Tata Institute of Fundamental Research, Homi Bhabha Road, Mumbai~400~005, India\\
\noindent $^{\ddagger}$UMR 7539 du CNRS, D\'epartement de Math\'ematiques, Universit\'e Paris 13, F-93430 Villetaneuse, France\\
\noindent Email: eva.bayer@epfl.ch; parimala@math.tifr.res.in;
queguin@math.univ-paris13.fr}

\volume{113}

\mon{November}

\parts{4}

\Date{MS received 02 July 2003}

\begin{abstract}
The question of the existence of an analogue, in the framework of
central simple algebras with involution, of the notion of Pfister form
is raised. In particular, algebras with orthogonal involution which
split as a tensor product of quaternion algebras with involution are
studied. It is proven that, up to degree 16, over any extension over
which the algebra splits, the involution is adjoint to a Pfister form.
Moreover, cohomological invariants of those algebras with involution are
discussed.
\end{abstract}


\maketitle

\setcounter{section}{-1}

\section{Introduction}

An involution on a central simple algebra is nothing but a twisted form
of a symmetric or alternating bilinear form up to a scalar factor
(\cite{KMRT}, ch.~1). Hence the theory of central simple algebras with
involution naturally appears as an extension of the theory of quadratic
forms, which is an important source of inspiration for this subject.

We do not have, for algebras with involution, such a nice algebraic
theory as for quadratic forms, since orthogonal sums are not always
defined, and are not unique when defined~\cite{D}. Nevertheless, in view
of the fundamental role played by Pfister forms in the theory of
quadratic forms, and also of the nice properties they share, it seems
natural to try and find out whether an analogous notion exists in the
setting of algebras with involution. 

The main purpose of this paper is to raise this question, which was
originally posed by David Tao~\cite{Tao:Pfister}; this is done in
\S2. In particular, this leads to the
consideration of algebras with orthogonal involution which split as a
tensor product of $r$ quaternion algebras with involution. One central
question is then the following: consider such a product of quaternions
with involution, and assume the algebra is split. Is the corresponding
involution adjoint to a Pfister form? The answer is positive up to
$r=5$. A survey of this question is given in \S2.4.
In \S4, we give a direct proof of this fact for
$r=4$. Before that, we study in \S3, the
existence of cohomological invariants for some of the algebras with
involution which can naturally be considered as generalisations of
Pfister quadratic forms.

\section{Notations}\label{not.section} 

Throughout this paper, the base field $F$ is supposed to be 
of characteristic different from $2$, and $K$ denotes a field 
extension of $F$. We refer the reader to~\cite{Sch}, \cite{Lam} and
\cite{KMRT} for more details on what follows in this section. 

\subsection{\it Cohomology} 

Let $K_s$ be a separable closure of the field $K$, and let us denote by
$\Gamma_K$ the absolute Galois group $\Gamma_K={\rm Gal}(K_s/K)$. The
Galois cohomology groups of $\Gamma_K$ with coefficients in $\mz/2$ will
be denoted by $H^i(K)=H^i(\Gamma_K,\mz/2)$. For any $a \in
K^\star$, we denote by $(a)$ the image in $H^1(K)$ of the class of
$a$ in $K^\star/{K^{\star 2}}$ under the canonical isomorphism
$K^\star/{K^{\star 2}}\simeq H^1(K)$, and by $(a_1,\dots,a_i)\in H^i(K)$
the cup-product $(a_1)\cup(a_2)\cup \dots\cup(a_i)$. In particular, the
element $(a_1,a_2)\in H^2(K)$ corresponds, under the canonical
isomorphism $H^2(K)\simeq {{\br}_2}(K)$, where ${{\br}_2}(K)$ denotes the
$2$-torsion part of the Brauer group of $K$, to the Brauer class of the
quaternion algebra $(a_1,a_2)_K$. 

Consider now a smooth integral variety $X$ over $F$, and denote by
$F(X)$ its function field. An element $\alpha\in H^i(F(X))$ is said to
be {\em unramified} if for each codimension one point $x$ in $X$, with
local ring ${\mathcal O}_x$ and residue field $\kappa_x$, the element
$\alpha$ belongs to the image of the natural map $H^i_{et}({\mathcal
O}_x,\mz/2)\ra H^i(F(X))$, or equivalently its image under the residue
map $\partial_x:\,H^i(F(X))\ra H^{i-1}(\kappa_x)$ is zero (see 
\cite{CT}, Theorem~4.1.1). We denote by $H^i_{nr}(F(X)/F)$ the subgroup
of $H^i(F(X))$ of unramified elements.

\subsection{\it Quadratic forms}\label{quadform.section}

The quadratic forms considered in this paper are non-degenerate. 
If $q$ is a quadratic form over $K$, we let $K(q)$ be the function field 
of the corresponding projective quadric. The field $K(q)$ is the generic 
field over which an anisotropic form $q$ acquires a non-trivial\break zero. 

Consider a diagonalisation $\qform{a_1,\dots,a_n}$ of a quadratic form
$q$. We denote by $d(q)$ the signed discriminant of $q$, that is
\hbox{$d(q)=(-1)^\frac{n(n-1)}{2}a_1\dots a_n\in\sq K$}, and by $C(q)$ its
Clifford algebra (see~\cite{Sch} or~\cite{Lam} for a definition and
structure theorems). We recall that $C(q)$ is a $\mz/2$-graded algebra,
and we denote by $C_0(q)$ its even\break part.

For any \hbox{$a_1,\dots,a_r\in K^\star$} we denote by $\pform{a_1,\dots,a_r}$
the $r$-fold Pfister form $\ot_{i=1}^r\qform{1,-a_i}$. We let $P_r(K)$
be the set of $r$-fold Pfister forms over $K$, and $GP_r(K)$ be the set
of quadratic forms over $K$ which are similar to an $r$-fold Pfister
form. Pfister forms are also characterized, up to similarities, by the
following property:

\begin{theor}[\!]\hskip -.65pc{\rm ((\cite{Kne}, Theorem~5.8) and \cite{W})}.\ \ 
Let $q$ be a quadratic form over $F$. The following assertions are
equivalent{\rm :}
\begin{enumerate}
\renewcommand{\labelenumi}{\rm (\roman{enumi})}
\leftskip .1pc
\item The dimension of $q$ is a power of $2${\rm ,} and $q_{F(q)}$ is
hyperbolic{\rm ;}

\item The quadratic form $q$ is similar to a Pfister form.
\end{enumerate}
\end{theor}

\pagebreak

From the above theorem, one easily deduces:

\begin{coro}\label{isoimphyp.cor}$\left.\right.$\vspace{.5pc}

\noindent Let $q$ be a quadratic form over $F$. The following assertions
are equivalent{\rm :}
\begin{enumerate}
\renewcommand{\labelenumi}{\rm (\roman{enumi})}
\leftskip .1pc
\item The dimension of $q$ is a power of\, $2$ and for any field extension
$K/F${\rm ,} if $q_K$ is isotropic{\rm ,} then it is hyperbolic{\rm ;}

\item $q$ is similar to a Pfister form.
\end{enumerate}\vspace{-.5pc}
\end{coro} 

We denote by $e_r$ the map \hbox{$P_r(K)\ra H^r(K)$} defined by Arason~\cite{A} 
as $e_r(\langle\langle a_1,\dots,$ $a_r\rangle\rangle)=(a_1,\dots ,a_r)$. 

Let $W(K)$ be the Witt ring of the field $K$, and denote by $I(K)$ the
fundamental ideal of $W(K)$, which consists of classes of
even-dimensional quadratic forms. Its $r$th power $I^r(K)$ is
additively generated by $r$-fold Pfister forms. For $r=1,2$ and $3$, the
invariant $e_r$ extends to a surjective homomorphism $I^r(K)\ra H^r(K)$
with kernel $I^{r+1}(K)$ (see~\cite{Mer:81} for $r=2$ and~\cite{MS} for
$r=3$). It follows from this that the class of an even-dimensional
quadratic form $q$ belongs to $I^2(K)$ (resp. $I^3(K)$) if and only if
$e_1(q)=0$ (resp. $e_1(q) = 0,$ $e_2(q)=0$). If we assume moreover that $q$ is
of dimension $4$ (resp. $8$), this is equivalent to saying that $q$ is
similar to a Pfister form.

Moreover, the maps $e_1$ and $e_2$ are actually defined (as maps) over
the whole Witt ring $W(K)$, and can be explicitly described in terms of
classical invariants of quadratic forms. Indeed, $e_1$ associates to the
class of a quadratic form $q$ its signed discriminant $d(q)\in\sq K$.
Moreover, the image under $e_2$ of the class of the same form $q$ is the
Brauer class of its Clifford algebra $C(q)$ if the dimension of $q$ is
even, and the Brauer class of $C_0(q)$ if the dimension of $q$ is odd. 

\subsection{\it Algebras with involution}\label{algainv.section}

An involution $\tau$ on a central simple algebra $B$ over $K$ is an
anti-automorphism of order $2$ of the ring $B$. We only consider here
involutions of the first kind, that is $K$-linear ones. For any field
extension $L/K$, we denote by $B_L$ the $L$-algebra $B\otimes_K L$, by
$\tau_L$ the involution $\tau\ot{\id}$ of $B_L$ and by $(B,\tau)_L$ the
pair $(B_L,\tau_L)$. 

Consider now a splitting field $L$ of $B$, that is an extension $L/K$
such that $B_L$ is the endomorphism algebra of some $L$-vector space
$V$. The involution $\tau_L$ is the adjoint involution ${\ad}_b$ with
respect to some bilinear form $b:\,V\times V\ra L$, which is either
symmetric or skew-symmetric. The type of the form $b$ does not depend on
the choice of the splitting field $L$; the involution $\tau$ is said to
be of orthogonal type if $b$ is symmetric, and of symplectic type if it
is skew-symmetric. 

Let $Q$ be a quaternion algebra over $K$. It admits a unique involution
of symplectic type, which we call the canonical involution of $Q$, and
which is defined by $\gamma_Q(x)={\trd}_Q(x)-x$, where ${\trd}_Q$ is the
reduced trace on $Q$. We denote by $Q^0$ the subspace of pure
quaternions, that is those $q\in Q$ satisfying ${\trd}_Q(q)=0$, or
equivalently, $\gamma_Q(q)=-q$. For any pure quaternion $q\in Q^0$, we
have $q^2\in F$. For any orthogonal involution $\sigma$ on $Q$, there
exists a pure quaternion $q\in Q$ such that
$\sigma={\Int}(q)\circ\gamma_Q$, where ${\Int}(q)$ is the inner automorphism
associated to $q$, defined by ${\Int}(q)(x)=qxq^{-1}$.

If the degree of $B$ is even, and if $\tau$ is of orthogonal type, we
denote by $d(\tau)\in\sq K$ the discriminant of $\tau$, and by
$C(B,\tau)$ its Clifford algebra~(\cite{KMRT}, \S7, 8). In the split
orthogonal case $(B,\tau)=({\End}_F(V),{\ad}_q)$, they correspond
respectively to the discriminant of $q$ and its even Clifford algebra
$C_0(q)$. Note that by the structure theorem~(\cite{KMRT}, (8.10)), if the
discriminant of $\tau$ is trivial, then the Clifford algebra $C(B,\tau)$
is a direct product of two $K$-central simple algebras,
$C(B,\tau)=C_+\times C_-$. 

A right ideal $I$ of a central simple algebra with involution $(B,\tau)$
is called isotropic if $\sigma(I)I=\{0\}$. The algebra with involution
$(B,\tau)$ is called isotropic if it contains a non trivial isotropic
right ideal, and hyperbolic if it contains a non trivial isotropic right
ideal of maximal dimension (that is of reduced dimension $\frac 12
\deg(B)$)~(\cite{KMRT}, \S~6).

In \cite{Tao:94}, David Tao associates to an algebra with orthogonal
involution $(B,\tau)$ a variety which, in the split orthogonal case
$(B,\tau)=({\End}_K(V),{\ad}_q)$, is the projective quadric associated to
$q$. This variety is called the involution variety of $(B,\tau)$; its
function field is the generic field over which $B$ splits and $\tau$
becomes isotropic.

\section{Three classes of algebras with involution}\label{classes.section}

From now on, we consider a central simple algebra $A$ over $F$, endowed
with an involution $\sigma$ of orthogonal type. We denote by $F_A$ the
function field of the Severi--Brauer variety of $A$, which is known to be
a generic splitting field for $A$. After scalar extension to $F_A$, the
involution $\sigma$ becomes the adjoint involution with respect to some
quadratic form over $F_A$, which we denote by $q_\sigma$. Note that this
form is uniquely defined up to a scalar factor in $F_A^\star$.

In view of the definition and properties of Pfister forms, it seems
natural, for our purpose, to consider the three classes of algebras with
involution introduced in this section.

\subsection{\it Pfister involutions} 

\setcounter{theo}{0}
\begin{defn}$\left.\right.$\vspace{.5pc}

\noindent {\rm The algebra with orthogonal involution $(A,\sigma)$ is called a Pfister
algebra with involution if $\sigma_{F_A}$ is adjoint to a Pfister form.}
\end{defn}

\begin{rema}\label{pfiinv.remark}
$\left.\right.$
{\rm \begin{enumerate}
\renewcommand{\labelenumi}{(\roman{enumi})}
\leftskip .5pc
\item If $(A,\sigma)$ is a Pfister algebra with involution, 
the degree of $A$ is a power of $2$. 

\item Since the form $q_\sigma$ is uniquely defined up to a scalar
factor, $(A,\sigma)$ is a Pfister algebra with involution if and only if
$q_\sigma\in GP(F_A)$. Moreover, two similar Pfister forms are actually
isometric (as follows from~(\cite{Sch}, ch.~4, 1.5). Hence, in this
particular case, there is a canonical choice for the quadratic form
$q_\sigma$; we may assume it is a Pfister form, in which case it is
uniquely defined up to isomorphism. 

\item Since any $2$-dimensional quadratic form is similar to a Pfister form, 
any degree $2$ algebra with orthogonal involution is a Pfister algebra 
with involution.
\end{enumerate}}
\end{rema}

The field $F_A$ is a generic splitting field for $A$. Hence, we may
deduce from the definition and from Corollary~\ref{isoimphyp.cor} the
following proposition:

\begin{propo}\label{pfiinv.prop}$\left.\right.$\vspace{.5pc}

\noindent The following assertions are equivalent{\rm :}
\begin{enumerate}
\renewcommand{\labelenumi}{\rm (\roman{enumi})}
\leftskip .5pc
\item $(A,\sigma)$ is a Pfister algebra with involution{\rm ;}

\item For any field extension $K/F$ which splits $A${\rm ,} the involution 
$\sigma_K$ is adjoint to a Pfister form{\rm ;}

\item The degree of $A$ is a power of\, $2$ and 
for any field extension $K/F$ which splits $A${\rm ,} if $\sigma_K$ 
is isotropic{\rm ,} then it is hyperbolic{\rm ;}

\item The degree of $A$ is a power of\, $2$ and 
after extending scalars to the function field of its involution 
variety{\rm ,} $(A,\sigma)$ becomes hyperbolic.
\end{enumerate}
\end{propo}

\subsection{\it Involutions of type $I\Rightarrow H$} 

As recalled in \S1, `isotropy implies hyperbolicity' is a
characterization of Pfister forms. Hence, we may also consider algebras
with involution satisfying the same property: 

\begin{defn}$\left.\right.$\vspace{.5pc}

\noindent {\rm The algebra with orthogonal involution $(A,\sigma)$ is
said to be of type $I\Rightarrow H$ if the degree of $A$ is a power of
$2$ and for any field extension $K/F$, if $(A,\sigma)_K$ is isotropic,
then it is hyperbolic.}
\end{defn} 

\begin{rema}\label{iimph.rem}$\left.\right.$
{\rm \begin{enumerate}
\renewcommand{\labelenumi}{(\roman{enumi})}
\leftskip .1pc
\item Again the condition is empty in degree $2$. Any degree $2$ algebra 
with orthogonal involution is of type $I\Rightarrow H$. 

\item From the previous proposition, one deduces that any involution 
of type $I\Rightarrow H$ is a Pfister involution. 
Moreover, if $A$ is split, then the two definitions 
are equivalent.
\end{enumerate}}
\end{rema} 

\subsection{\it Product of quaternions with involution} 

Up to similarities, Pfister forms are those quadratic forms which
diagonalise as a tensor product of two dimensional forms. Hence, we
now consider algebras with involution which split as a tensor product of
degree $2$ algebras with involution.

\begin{defn}$\left.\right.$\vspace{.5pc}

\noindent {\rm The algebra with orthogonal involution $(A,\sigma)$ is called a product
of quaternions with involution if there exists an integer $r$ and
quaternion algebras with involution $(Q_i,\sigma_i)$ for $i=1,\dots, r$
such that $(A,\sigma)\simeq \otimes_{i=1}^r(Q_i,\sigma_i)$.}
\end{defn} 

\begin{rema}\label{proqua.rem}$\left.\right.$
{\rm \begin{enumerate}
\renewcommand{\labelenumi}{(\roman{enumi})}
\leftskip .5pc
\item If $(A,\sigma)$ is a product of quaternions with involution, then
the degree of $A$ is a power of $2$. 

\item Since $\sigma$ is of orthogonal type, the number of indices $i$
for which $\sigma_i$ is of symplectic type is necessarily even. 

\item In \cite{KPS:91}, it is proven that a tensor product of two
quaternion algebras with orthogonal involutions admits a decomposition
as a tensor product of quaternion algebras with symplectic (hence
canonical) involutions. Hence, any product of quaternions with
involution admits a decomposition as above in which all the $\sigma_i$
if $r$ is even, and all but one if $r$ is odd, are the canonical
involutions of $Q_i$.
\end{enumerate}}
\end{rema} 

As above, the condition is empty in degree $2$, any degree $2$ algebra
with orthogonal involution is a product of quaternions with involution.
In degrees $4$ and $8$, we have the following characterizations:

\begin{theor}[\!]\hskip -.65pc{\rm \cite{KPS:91}.}\ \ \label{KPS.thm}
Let $(A,\sigma)$ be a degree $4$ algebra with orthogonal 
involution. It is a product of quaternions with involution if and only if 
the discriminant of $\sigma$ is $1$. 
\end{theor}

\begin{theor}[\!]\hskip -.65pc{\rm (\cite{KMRT}, (42.11)).}\ \ \label{KPS8.thm}
Let $(A,\sigma)$ be a degree $8$ algebra with orthogonal involution. It
is a product of quaternions with involution if and only if the
discriminant of $\sigma$ is trivial and one component of the Clifford
algebra of $(A,\sigma)$ splits.
\end{theor}

\subsection{\it Shapiro's conjecture}\label{sha.section}

It seems a natural question to try and find out whether the three
classes of algebras with involution introduced above are equivalent.
This is obviously the case in degree $2$. The following proposition will
be proven in \S3.3:

\begin{propo}\label{deg8.prop}$\left.\right.$\vspace{.5pc}

\noindent Let $(A,\sigma)$ be an algebra of degree at most $8$ 
with orthogonal involution. The following are equivalent{\rm :}
\begin{enumerate}
\renewcommand{\labelenumi}{\rm (\roman{enumi})}
\leftskip .5pc
\item $(A,\sigma)$ is a Pfister algebra with involution{\rm ;}

\item $(A,\sigma)$ is a product of quaternions with involution{\rm ;}

\item $(A,\sigma)$ is of type $I\Rightarrow H$.
\end{enumerate}
\end{propo}

Nevertheless, the general question of the equivalence of these three
classes of algebras with involution is largely open in higher degree.
The most significant result is due to Shapiro. In his book
`Composition of quadratic forms' he makes the following conjecture:

\begin{conjecture}\hskip -.4pc{\rm (\cite{Sha}, (9.17)).\ \ \label{sha.conj} 
Let $(A,\sigma)$ be a product of $r$ quaternions with involution. 
If $A$ is split, then $(A,\sigma)$ admits a decomposition 
as a tensor product of $r$ quaternion algebras with involution 
in which each quaternion algebra is split.}
\end{conjecture}

Moreover, he proves the following theorem:

\begin{theor}[\!]\hskip -.5pc{\rm (\cite{Sha}, Claim in p.~166 and Ch.~9).}\ \ 
Conjecture~{\rm \ref{sha.conj}} is true if $r\leq5$.
\end{theor} 

It is easy to see that Shapiro's conjecture is true for some $r$ if and
only if any product of $r$ quaternions with involution is a Pfister
algebra with involution. Hence the previous theorem implies.

\begin{coro}$\left.\right.$\vspace{.5pc}

\noindent Any product of $r\leq 5$ quaternions with involution is a
Pfister algebra with involution. 
\end{coro}

Shapiro does not give a direct proof of this conjecture. He is actually
interested in another conjecture, which he calls the Pfister factor
conjecture, and which gives a characterization of $r$-fold Pfister forms
in terms of the existence of vector-spaces of maximal dimension in the
group of similarities of these forms (see~\cite{Sha}, (2.17)) for a
precise statement). He proves the Pfister factor conjecture for $r\leq
5$, using tools from the algebraic theory of quadratic forms, and also
proves it is equivalent to Conjecture~\ref{sha.conj}.

In fact, for $r\leq 3$, we have a little bit more: as already mentioned
in Proposition~\ref{deg8.prop}, $(A,\sigma)$ is a product of quaternions
with involution if and only if it is a Pfister algebra with involution.
A proof of this fact, using cohomological invariants, which was already
noticed by David Tao, will be given in \S3,
where we study the general question of cohomological invariants of
Pfister involutions.

In \S4, we give a direct proof of~\ref{sha.conj} in the
$r=4$ case, based on the study of some trace forms of product of
quaternions with involution. Since this paper was submitted, Serhir
and Tignol~\cite{ST} found another direct proof of this conjecture
for $r\leq 5$, using the discriminant of symplectic involutions defined
by Berhuy, Monsurro and Tignol~\cite{BMT}.

\section{Cohomological invariants}\label{cohinv.section} 

From the point of view of quadratic form theory, cohomological
invariants seem a natural tool for studying these questions. In
\S~\ref{ei.section}, we define an invariant of a Pfister involution,
with values in the unramified cohomology group of the function field of
the generic splitting field of the underlying algebra. We then study the
question of the existence of an analogous invariant with values in the
cohomology group of the base field. 

\subsection{\it Invariant $e_i$ for Pfister algebras with involution}\label{ei.section} 

Throughout this section, $(A,\sigma)$ is a Pfister algebra with
involution over $F$. The degree of $A$ is $2^i$, and we assume
$q_\sigma$ is an $i$-fold Pfister form over $F_A$ (see
Remark~\ref{pfiinv.remark}). Let us consider the Arason invariant
$e_i(q_\sigma)\in H^i(F_A)$. We have the following:

\setcounter{theo}{0}
\begin{theor}[\!]\label{unr.thm}
The invariant $e_i(q_\sigma)$ belongs to the unramified cohomology 
group $H^i_{nr}(F_A/F)$. 
\end{theor}

\begin{proof}
Given a codimension one point $x$ of the Severi--Brauer variety $X_A$ of
$A$, its residue field $\kappa_x$ splits $A$. Hence, the involution
$\sigma_{\kappa_x}$ is the adjoint involution with respect to a
quadratic form $q_x$ which is a Pfister form uniquely determined by
$\sigma_{\kappa_x}$ (see Remark~\ref{pfiinv.remark}~(ii)).

Let us now consider the completions $\widehat{{\mathcal O}_x}$ and
$\widehat{F(X_A)}$ of ${\mathcal O}_x$ and $F(X_A)$ at the discrete
valuation associated to $x$. Since $\widehat{{\mathcal O}_x}$ is
complete, the field $\widehat{F(X_A)}$ is isomorphic to $\kappa_x((t))$,
and for the same reason as above, the involution
$\sigma_{\widehat{F(X_A)}}$ is adjoint to a unique Pfister form
$q_{\widehat{F(X_A)}}$, which is the form $q_x$ extended to
$\kappa_x((t))$.

From this, we get that $e_i(q_{\widehat{F(X_A)}})$ is the image of
$e_i(q_x)$ under the natural map $H^i(\kappa_x)\ra H^i(\kappa_x((t)))$.
By~(\cite{CT}, \S3.3), since the corresponding ring is complete, this
implies that the image of $e_i(q_{\widehat{F(X_A)}})$ under the residue
map \hbox{$\partial_x:\,H^i(\widehat{F(X_A)})\ra H^{i-1}(\kappa_x)$} is
trivial. Finally, again by (\cite{CT}, \S3.3),
$\partial_x(e_i(q_\sigma))=\partial_x(e_i(q_{\widehat{F(X_A)}}))$, and
this proves the theorem.

\hfill \ab
\end{proof} 

Of course, it would be nicer to have an invariant with values in the
cohomology group of the base field. To be more precise, let us denote by
$E_i(A)$ the kernel of the restriction map $H^i(F)\ra H^i(F_A)$ and by
$\Phi$ the injection:
\begin{equation*}
\Phi:\,H^i(F)/E_i(A)\ra H^i_{nr}(F_A).
\end{equation*}
We may ask the following question: Does $e_i(q_\sigma)$ belong to the
image of $\Phi$? In \S\ref{i=012.section}, we prove that this is
the case for $i=0,1$ and $2$, and we give an interpretation of the
corresponding invariant in $H^i(F)/E_i(A)$ in terms of classical
invariants of orthogonal involutions. In \S\ref{noe3.section}, we
prove this is not the case anymore for $i=3$.

\subsection{\it Invariants $e_0$, $e_1$ and $e_2$}\label{i=012.section}

Let us consider now any algebra with orthogonal involution $(A,\sigma)$.
As recalled in \S1, the first three Arason
invariants $e_0, e_1$ and $e_2$ for quadratic forms play a particular
role. Indeed, they are actually defined as maps over the whole Witt ring
$W(F)$, and they can be described in terms of classical invariants of
quadratic forms. In view of this, we may give the following definition:

\begin{defn}$\left.\right.$\vspace{.5pc}

\noindent {\rm Let $(A,\sigma)$ be an algebra with orthogonal involution
over $F$. We let 
\begin{equation*}
e_0(A,\sigma)=\overline{\deg(A)}\in \mz/2\mz\simeq H^0(F).
\end{equation*}

If the degree of $A$ is even (that is $e_0(A,\sigma)=0$), 
we let 
\begin{equation*}
e_1(A,\sigma)=d(\sigma)\in\sq F=H^1(F),
\end{equation*}
where $d(\sigma)$ denotes the discriminant of $\sigma$.}
\end{defn}

\begin{rema}{\rm 
Note that, as opposed to what happens for quadratic forms, the invariant
$e_1$ is only defined when $e_0$ is trivial. This is a consequence of
the fact that the discriminant of a quadratic form is an invariant up to
similarity, and hence an invariant of the corresponding adjoint
involution, only if the form has even dimension.}
\end{rema} 

Assume now that $e_0(A,\sigma)=e_1(A,\sigma)=0$, which means
$(A,\sigma)$ has even degree and trivial discriminant. From the
structure theorem recalled in \S1, the Clifford
algebra $C(A,\sigma)$ is isomorphic to a direct product of two central
simple algebras over $F$, $C(A,\sigma)=C_+\times C_-$, which give rise
to two Brauer classes $[C_+]$ and $[C_-]$ in ${\br}_2(F)$. The definition
of $e_2$ then relies on the following proposition:

\begin{propo}{\rm (\cite{KMRT}, (9.12))}$\left.\right.$\vspace{.5pc}

\noindent In ${\br}_2(F)$, we have $[C_+]+[C_-]\in\{0,[A]\}$. 
\end{propo}

Indeed this implies that the two classes actually coincide in the
quotient of ${\br}_2(F)$ by the subgroup $\{0,[A]\}$, which is exactly
$E_2(A)$. Hence, we give the following definition:

\begin{defn}$\left.\right.$\vspace{.5pc}

\noindent {\rm Let $(A,\sigma)$ be an algebra with orthogonal involution over $F$ of
even degree and trivial discriminant. We let
\begin{equation*}
e_2(A,\sigma)=[C_+]=[C_-]\in{\br}_2(F)/E_2(A).
\end{equation*}}
\end{defn}

Next, we prove the following:

\begin{propo}\label{e_i.prop}$\left.\right.$\vspace{.5pc}

\noindent Let $(A,\sigma)$ be a split algebra with orthogonal
involution{\rm ,} $(A,\sigma)=({\End}_F(V),{\ad}_q)$. When they are defined{\rm ,} the
invariants $e_0(A,\sigma), e_1(A,\sigma)$ and $e_2(A,\sigma)$ coincide
respectively with $e_0(q), e_1(q)$ and $e_2(q)$.
\end{propo}

\begin{proof}
This is clear for $e_0$ and $e_1$. For $e_2$, first note that if $A$ is
split, then $e_2(A,\sigma)$ actually belongs to ${\br}_2(F)$. Moreover,
$e_2(A,\sigma)$ is only defined when $e_0$ and $e_1$ are trivial, in
which case the form $q$ is of even dimension and trivial discriminant.
From the structure theorem for Clifford algebra (see for
instance~(\cite{Lam}, 5, \S2) or (\cite{Sch}, 9(2.10)) we get that in this
situation, we may represent $C(q)$ as $M_2(B)$, for some central simple
algebra $B$ over $F$, and the even part $C_0(q)$ corresponds to diagonal
matrices, $C_0(q)\simeq B\times B$, so that
$e_2(q)=[C(q)]=[B]=e_2(A,\sigma)$.\hfill \ab
\end{proof}

From this proposition, we easily deduce:

\begin{coro}$\left.\right.$\vspace{.5pc}

\noindent Let $(A,\sigma)$ be an algebra with orthogonal involution 
such that $e_i(A,\sigma)$ is defined for some $i\leq 2$. 
The invariant $e_i(A,\sigma)$ maps to $e_i(q_\sigma)$ 
under the morphism $\Phi:\,H^i(F)/E_i(A)\ra H^i(F_A)$.
\end{coro}

Hence, those invariants $e_i$ may be used to characterize degree $4$ and
$8$ Pfister involutions. Indeed, consider an algebra with orthogonal
involution $(A,\sigma)$, of degree $2^i$ for some $i\in\{2,3\}$. By
definition, it is a Pfister algebra with involution if and only if the
form $q_\sigma$ belongs to $GP_i(F_A)$. As recalled in
\S1, this is also equivalent to saying that
$e_1(q_\sigma)=0$ if $i=2$, and $e_1(q_\sigma)=e_2(q_\sigma)=0$ if
$i=3$. From this we get the following:

\begin{propo}\label{above.prop}$\left.\right.$\vspace{.5pc}

\noindent The degree $4$ algebra with orthogonal involution $(A,\sigma)$
is a Pfister algebra with involution if and only if $e_1(A,\sigma)=0$.
The degree $8$ algebra with orthogonal involution $(A,\sigma)$ is a
Pfister algebra with involution if and only if
$e_1(A,\sigma)=e_2(A,\sigma)=0$.
\end{propo} 

Using this, we are now able to prove Proposition~\ref{deg8.prop}.

\subsection{\it Proof of Proposition~\ref{deg8.prop}}\label{proof.section}

Comparing Proposition~\ref{above.prop} with Theorems~\ref{KPS.thm}
and \ref{KPS8.thm} we get the equivalence between (i) and (ii), using
(\cite{KMRT}, (9.14)) in the degree 8 case. Moreover, as already noticed
in Remark~\ref{iimph.rem}, any involution of type $I\Rightarrow H$ is a
Pfister involution. Hence, it only remains to prove that a product of
$r$ quaternions with involution with $r\leq 3$ is of type $I\Rightarrow
H$. Let us consider such a product of quaternions with involution
$(A,\sigma)$ and assume it is isotropic. Then, $A$ cannot be a division
algebra and has index at most $2^{r-1}$.

If $A$ is split, $\sigma$ is the adjoint involution with respect to an
isotropic Pfister form. Hence it is hyperbolic, and this concludes the
proof in that case.

Assume now that the index of $A$ is $2^{r-1}$, and let $D$ be a division
algebra Brauer-equivalent to $A$. We may represent $(A,\sigma)$ as
$({\End}_{D}(M),{{\ad}}_h)$, where $(M,h)$ is a rank $2$ hermitian module over
$D$. Again, since $\sigma$ is isotropic, $h$ is isotropic, hence
hyperbolic because of its rank, and this concludes the proof in that
case. If $r=2$, we are done, and it only remains to consider the case when
$r=3$ and $A$ has index $2$. Let $Q$ be a quaternion division algebra
Brauer-equivalent to $A$, denote by $\gamma$ its canonical involution,
and let $(M,h)$ be a skew-hermitian module over $(Q,\gamma)$ such that
$(A,\sigma)=({\End}_{Q}(M),{\ad}_h)$. Denote by $C$ the conic associated to
$Q$, and by $L$ its function field, which is known to be a generic
splitting field for $Q$, and hence for $A$. Since $A_L$ is split,
$\sigma_L$, and hence $h_L$ are hyperbolic. By (\cite{PSS}, Proposition~3.3)
(see also \cite{Dej:01}), this implies that $h$ itself, and hence
$\sigma$ is hyperbolic, and the proof is complete.

\subsection{\it About the $e_3$ invariant}\label{noe3.section}

As opposed to what happens for $e_0$, $e_1$ and $e_2$, there does not
exist any invariant in $H^3(F)/E_3(A)$ which is a descent of
$e_3(q_\sigma)$ for degree 8 Pfister algebras with involution, as
the following theorem shows:

\begin{theor}[\!]\label{noe3.thm} 
There exists a degree $8$ Pfister algebra with involution for which the
invariant $e_3(q_\sigma)$ does not belong to the image of the morphism
$\Phi:\,H^3(F)/E_3(A)\ra H^3_{nr}(F(X_A))$.
\end{theor}

\begin{proof} 
In his paper `Simple algebras and quadratic forms', Merkurjev
(\cite{Mer:92}, proof of Theorem~4) constructs a division algebra $A$,
which is a product of three quaternion algebras, $A=Q_1\otimes
Q_2\otimes Q_3$, and with centre a field $F$ of cohomological dimension
at most $2$. In particular, we have $H^3(F)=0$. Consider any orthogonal
decomposable involution $\sigma=\sigma_1\otimes\sigma_2\otimes\sigma_3$
on $A$. By Proposition~\ref{deg8.prop}, $(A,\sigma)$ is a degree $8$
Pfister algebra with involution. Moreover, by a result of
Karpenko~(\cite{Kar:00}, Theorem~5.3), since $A$ is a division algebra, the
involution $\sigma$ remains anisotropic over $F(X_A)$. Hence, $q_\sigma$
is an anisotropic 3-fold Pfister form, and $e_3(q_\sigma)$ is 
non-trivial. Since $H^3(F)=0$, this is enough to prove that $e_3(q_\sigma)$
does not belong to the image of $\Phi$.\hfill \ab
\end{proof} 

\begin{rema}{\rm 
Using Merkurjev's construction of division product of quaternions with
involution mentioned in the proof of Theorem~\ref{noe3.thm}, one may
construct explicit elements in the unramified cohomology
$H^i_{nr}(F_A/F)$ for any $i\geq 3$ for which Shapiro's conjecture is
known, which do not come from $H^i(F)$.}
\end{rema}

\section{Product of four quaternions with involution}\label{r=4.section}

In this section, we give a direct proof of Shapiro's conjecture for
$r=4$, i.e. we prove that any product of four quaternions with
involution is a Pfister algebra with involution.

By Proposition~\ref{pfiinv.prop} and Corollary~\ref{isoimphyp.cor}, it
suffices to prove the following proposition:

\setcounter{theo}{0}
\begin{propo}\label{r=4.prop}$\left.\right.$\vspace{.5pc}

\noindent Let $(A,\sigma)$ be a product of four quaternions with
involution. If $A$ is split and $\sigma$ is isotropic{\rm ,} then it is
adjoint to a hyperbolic quadratic form.
\end{propo} 

Let $(A,\sigma)=\otimes_{i=1}^4(Q_i,\gamma_i)$, and assume $A$ is split
and $\sigma$ is isotropic. By Remark~\ref{proqua.rem} (ii), we may
assume that each $\gamma_i$ is the canonical involution on $Q_i$. Let us
denote by $(D,\gamma)=(Q_1,\gamma_1)\otimes (Q_2,\gamma_2)$. We start
with a lemma which gives a description of $(A,\sigma)$:

\begin{lem}\label{iso.lem}
There exists an invertible element $u\in D^\star$ satisfying
$\gamma(u)=u$, ${\trd}_D(u)=0$ and ${\nrd}_D(u)\in F^{\star 2}$ such that
\begin{equation*}
(Q_3,\gamma_3)\otimes (Q_4,\gamma_4)\simeq (D,{\Int}(u^{-1})\circ\gamma)
\end{equation*}
and
\begin{equation*}
(A,\sigma)\simeq({\End}_F(D),{\ad}_{q_u}),
\end{equation*}
where $q_u$ is the quadratic form defined on $D$ by
$q_u(x)={\trd}_D(xu\gamma(x))$.
\end{lem} 

\begin{proof}
Since $A$ is split, $Q_3\otimes Q_4$ is isomorphic to $D$, and
$\gamma_3\otimes\gamma_4$ corresponds under this isomorphism to an
orthogonal involution $\gamma'$ on $D$. There exists an invertible
$\gamma$-symmetric element $u\in D$ such that $\gamma'=
{\Int}(u^{-1})\circ\gamma$. Moreover, since $\gamma_3\otimes\gamma_4$ is
decomposable, by Theorem~\ref{KPS.thm}, its discriminant is trivial.
Hence, so is the discriminant of $\gamma'$, and
by~(\cite{KMRT}, (7.3)(1)), we get that ${\nrd}_D(u)\in F^{\star 2}$.

Using this, we now get that $(A,\sigma)$ is isomorphic to 
\begin{equation*}
(D\otimes D,\gamma\otimes{\Int}(u^{-1})\circ\gamma).
\end{equation*}
By~(\cite{KMRT}, (11.1)), under the canonical isomorphism $D\otimes
D\simeq {\End}_F(D)$, the involution
$\gamma\otimes{\Int}(u^{-1})\circ\gamma$ is adjoint to the quadratic form
$q_u:\,D\ra F$ defined by $q_u(x)={\trd}_D(xu\gamma(x))$, and it only
remains to prove that we may assume ${\trd}_D(u)=0$.

Since the isotropic involution $\sigma$ is adjoint to $q_u$, the
quadratic form also is isotropic. Moreover, by a general position
argument, there exists an invertible element $y\in D$ such that
$q_u(y)={\trd}_D(yu\gamma(y))=0$. Then, the map $D\ra D$, $x\mapsto
xy^{-1}$ is an isometry between $q_u$ and the quadratic form
$q_{yu\gamma(y)}:\,x\mapsto{\trd}_D(xyu\gamma(y)\gamma(x))$. One may
easily check that this new element $yu\gamma(y)$ satisfies all
properties of the lemma, including ${\trd}_D(yu\gamma(y))=0$, and this
ends the proof.\hfill \ab
\end{proof} 

To get Proposition~\ref{r=4.prop}, we now have to prove that the
quadratic form $q_u$ is hyperbolic. This follows easily from the
following lemma:

\begin{lem}\label{isosub.lem}
The quadratic space $(D,q_u)$ contains a totally isotropic 
subspace of dimension $5$. 
\end{lem}

Indeed, by the computations of classical invariants for tensor product 
of algebras with involution given in~(\cite{KMRT}, (7.3)(4) and p.~150),
we have $e_0(A,\sigma)=e_1(A,\sigma)=e_2(A,\sigma)=0$. 
Hence, by Proposition~\ref{e_i.prop}, the quadratic form $q_u$ has trivial 
$e_0$, $e_1$ and $e_2$ invariants, and as recalled in 
\S1, this implies that it lies in $I^3(K)$. 

Since $q_u$ is $16$ dimensional, the previous lemma implies its
anisotropic dimension is at most $6$. By Arason--Pfister's theorem, this
implies that $q_u$ is hyperbolic, and thus concludes the proof of
Proposition~\ref{r=4.prop}.

\setcounter{theo}{2}
\begin{pol}{\rm 
For any $z\in D^\star$, we denote by $q_z$ the quadratic form 
$D=Q_1\otimes Q_2\ra F$, $x\mapsto{\trd}_D(xz\gamma(x))$. 
We first prove the following fact:

\begin{claim}\label{1.claim}  
{\rm Let $z\in D^\star$ satisfy ${\trd}_D(z)=0$. We then have
\begin{enumerate}
\renewcommand{\labelenumi}{(\roman{enumi})}
\leftskip .1pc
\item $Q_1$ is totally isotropic for $q_z$; 

\item for all $x\in Q_1$, ${\trd}_D(xz)=0$. 
\end{enumerate}}
\end{claim}

Indeed, for any $x\in Q_1$, we have $x\gamma(x)=x\gamma_1(x)={\nrd}_D(x)$.
Hence, $q_z(x)={\trd}_D(xzx^{-1}{\nrd}_D(x))={\nrd}_D(x){\trd}_D(z)=0$.
Moreover, considering the corresponding bilinear form, we also get that
for any $x,y\in Q_1$, ${\trd}_D(xz\gamma(y))=0$, and this gives the second
part of the claim by taking $y=1$.

Let us now denote by $\phi$ the isomorphism $(Q_3,\gamma_3)\otimes
(Q_4,\gamma_4)\simeq (D,{\Int}(u^{-1})\circ\gamma)$ of
Lemma~\ref{iso.lem}. For any pure quaternion $q\in Q_4^0$, we let $W_q$
be the image under $\phi$ of the $3$-dimensional subspace $\{x\otimes
q,\ x\in Q_3^0\}$ of $Q_3\otimes Q_4$. We then have

\begin{claim}\label{2.claim}
{\rm The subspace $\gamma(W_q)$ of $D$ is totally isotropic for $q_u$.}
\end{claim}

Indeed, from the corresponding properties for $\{x\otimes q,\ x\in
Q_3^0\}$, any element $y\in W_q$ satisfies $y^2\in F$ and
${\Int}(u^{-1})\circ\gamma(y)=y$. Hence, we have $\gamma(y)u=uy$ and
$q_u(\gamma(y)) ={\trd}_D(\gamma(y)uy)={\trd}_D(uy^2)=y^2{\trd}_D(u)=0$.

Let us now denote by $T$ the kernel of the linear form on $D$ defined by
$z\mapsto{\trd}_D(u\gamma(z))$. Clearly, $T\cap\gamma(W_q)$ has dimension
at least $2$. Fix a $2$ dimensional subspace $V_q\subset
T\cap\gamma(W_q)$. We then have

\begin{claim} 
{\rm The subspace $Q_1+V_q$ of $D$ is totally isotropic for $q_u$.}
\end{claim}

Indeed, $Q_1$ is totally isotropic by Claim~\ref{1.claim}, and since
$V_q\subset \gamma(W_q)$, it also is by Claim~\ref{2.claim}. Moreover,
for any $z\in V_q\subset T$, we have ${\trd}_D(u\gamma(z))=0$. Hence
by Claim~\ref{1.claim}(ii), ${\trd}_D(xu\gamma(z))=0$ for any $x\in Q_1$. Hence
$Q_1$ and $V_q$ are orthogonal, and we get the claim.

To finish with, it only remains to prove that there exists some $q\in
Q_4^0$ such that $Q_1+V_q$ has dimension greater than $5$, i.e. $V_q$ is
not contained in $Q_1$. But if $V_q$ is contained in $Q_1$, then it is
contained in $Q_1^0$ which has dimension $3$. One may then choose
another pure quaternion $q'\in Q_4^0$ which is linearly independent from
$q$. This way, we get another $2$ dimensional subspace $V_{q'}$ which is
in direct sum with $V_q$, and which cannot also be contained in
$Q_1^0$.}\hfill \ab
\end{pol}

\section*{Acknowledgement}

We would like to thank Jean-Pierre Tignol for his useful comments 
on a preliminary version of this paper.

\end{document}